\documentclass[11pt,leqno]{article}
\usepackage{amssymb,amsmath,amsfonts,amsbsy, xspace, latexsym, amscd, graphicx}

\newtheorem{theorem}{Theorem}[section]
\newtheorem{lemma}{Lemma}[section]

\newtheorem{definition}{Definition}[section]
\newtheorem{corollary}{Corollary}[section]

\newtheorem{simpler notation}{Simpler Notation}[section]
\newtheorem{remark}{Remark}[section]
\newtheorem{example}{Example}[section]

\newcommand{\Proof}{\noindent \textbf{Proof:\;}}

\newcommand{\eqdef}{\ensuremath{\overset{\text{def}}{=}}}

\begin{document}
{\small

\title{A Lost Theorem: Definite Integrals in An Asymptotic Setting}
\author{	Ray Cavalcante and Todor D. Todorov}
\date{}
\maketitle

\section{INTRODUCTION} \label{S: Introduction}

We present a simple yet rigorous theory of integration that is based on two axioms rather than on a construction involving Riemann sums. With several examples we demonstrate how to set up integrals in applications of calculus without using Riemann sums. In our axiomatic approach even the proof of the existence of  the definite integral (which does use Riemann sums) becomes slightly more elegant than the conventional one.  We also discuss an interesting connection between our approach and the history of calculus. The article is written for readers who teach calculus and its applications. It might be accessible to students under a teacher's supervision and suitable for senior projects on calculus, real analysis, or history of mathematics.

	Here is a summary of our approach. Let $\rho : [a,b] \to \mathbb{R}$ be a continuous function and let $I: [a, b] \times [a, b]\to\mathbb{R}$ be the corresponding  integral function, defined by
\[	
	I(x, y) = \int_{x}^{y}\rho(t)dt.
\]	
	Recall that $I(x, y)$ has the following two properties:
	\begin{description}
		\item{\textbf{(A)}} {\em Additivity}: $I(x, y) + I(y, z) = I(x, z)$ for all $x, y, z \in [a, b]$.
		\item{\textbf{(B)}} {\em Asymptotic Property}: $I(x, x + h) = \rho(x)h + o(h)$ as $h\to0$ for all $x \in [a, b]$, in the  sense that
		\[
			\lim_{h \to 0} \frac{I(x, x + h) - \rho(x)h}{h} = 0.
		\]
	\end{description}
	
	In this article we show that properties  \textbf{(A)} and \textbf{(B)} are characteristic of the definite integral. More precisely, we show that for a given continuous $\rho : [a, b] \to \mathbb{R}$, there is no more than one function $I: [a, b] \times [a, b] \to \mathbb{R}$ with properties \textbf{(A)} and \textbf{(B)}. This will justify the simple definition $\int_{a}^{b}\rho(x)dx\eqdef I(a, b) $, where $I(x, y)$ is a function satisfying \textbf{(A)} and \textbf{(B)}. In this manner, we are able to rigorously develop the theory of definite integrals and their applications without partitioning the interval $[a, b]$ and without using Riemann sums. Notice that, at this stage, the existence of the integral is still unsettled. Next, we prove that if $R(x)$ is an antiderivative of $\rho(x)$, then $I(a, b) = R(b) - R(a)$. Conventionally, this formula is used for explicit evaluation, but in our approach it also guarantees the existence of the integral for all $\rho(x)$ with an explicit antiderivative $R(x)$. Our approach utilizes the Riemann partition procedure for \emph{only one purpose}: to prove the existence of the definite integral for functions without explicit anti-derivatives such as $f(x)=e^{-x^{2}}$. Also, our proof (Theorem~\ref{T: General Existence Result}) seems to be slightly more elegant than the conventional one. Next we use axioms similar to {\bf(A)} and {\bf(B)} to define the concepts of area under the curve, arclength, volumes of revolution, etc. More precisely, we define  all of these geometrical quantities as being additive and asymptotically equal to their Euclidean  counterparts (such as the area of a rectangle, the Euclidean distance between two points, the volume of a cylindrical shell, etc.).  Our definitions are mathematically correct and well motivated. Derivations of the corresponding integral formulas (for the area under a curve, arclength, volume, etc.) appear in our approach as simple rigorous theorems with proofs done in the spirit of asymptotic analysis; we involve neither partitions of the interval $[a, b]$ nor Riemann sums.

	The elementary theory of integration presented in this article (and summarized above) opens the door for a simple yet completely rigorous method of teaching integration and its applications in a calculus course or a beginning real analysis course. Also, we strongly recommend our method of setting up integrals (without Riemann sums) for teaching  physics and engineering courses based on calculus.

	The method presented in this article has a long and interesting history.  The reader might be surprised to learn that practically all textbooks on calculus and its applications that were written in the period between Leibniz and Riemann motivate, define, and set up integrals by methods very similar to the method presented in this article, although disguised in the language of infinitesimals. The method (among other treasures) was lost in the history of calculus after infinitesimals were abolished at the end of 19th century.  This explains the choice of the title ``A Lost Theorem....'' We shall briefly discuss the connection of our approach with infinitesimal calculus in Section~\ref{S: Lost Theorem}. In the modern literature we identify three sources using methods for integration similar to our approach.  In H. J. Keisler's calculus textbook \cite{Keisler} (look for {\em Infinite Sum Theorem} on p. 303), the reader will find a method similar to ours in the framework of nonstandard analysis (in a very accessible form). We also refer to S. Lang~(\cite{Lang}, pp. 292-296),  L. Gillman and R. McDowell~(\cite{GillmanMc}, pp. 176-179) and  L. Gillman~\cite{Gillman}, where a  property similar to {\bf(B)} is used for a definition of the definite integral. 

\section{TOPICS IN ASYMPTOTIC ANALYSIS}\label{S: Topics in Asymptotic Analysis}
	For our axiomatic approach to integration we assume knowledge of limits, continuity, and derivatives at the level of a typical beginning calculus course. In contrast, we do not assume any knowledge of integration. Finally, we need several elementary concepts and notation borrowed from asymptotic analysis, which we present in this section. Most of the results are elementary, and we leave the proofs to the reader.
	\begin{definition} {\em We denote by $o(x^{n})$ the set of all functions $f: D_f\to\mathbb{R}$ such that $D_f\subseteq\mathbb{R},\; 0\in\overline{D_f}$, and $\lim_{x\to 0} f(x)/x^{n}=0$, where $\overline{D_f}$ stands for the closure of $D_f$ in $\mathbb{R}$. We summarize this as $o(x^{n}) = \{ f(x) : f(x)/x^{n} \to 0$ as x $\to 0\}$. It is customary to write $f(x) = o(x^{n})$ instead of the more precise $f \in o(x^{n})$ in the case when $f$ is an unspecified function in $o(x^n)$ . If $n = 0$, the above definition reduces to $f(x) = o(1)$ if $f(x) \to 0$ as $x \to 0$.
	}\end{definition}

	\begin{example} {\em
	$x^{2} = o(x)$ since $x^{2}/x \to 0$ as $x \to 0$. In contrast, $\sin x \ne o(x)$ since 
$\sin {x}/x \to 1$ as $x \to 0$. Also, we have $\sin {x} = o(1)$ since $\sin {x}\to 0$ as $x \to 0$. 
	}\end{example}
	\begin{theorem}[Increment Theorem] \label{T: Increment Theorem}   Let $f: [a, b]\to \mathbb{R}$ be a function and $x\in (a, b)$. Then $f$ is 
differentiable at $x$ if and only if $f(x+h)-f(x)=f^\prime(x)h+o(h)$.
		
	\end{theorem}
	\begin{lemma}	
		Let $f(x) = o(x^{m})$ and $g(x) = o(x^{n})$. Then {\bf (a)}  $f(x) \pm g(x) = o(x^{k})$ where $k = \min(m, n)$, and {\bf  (b)} $f(x)g(x) = o(x^{m+n})$. Also,   {\bf (c)}  if  $f(x) = o(1)$ and $g(x) = o(1)$ and $f$ is continuous at $0$, then $f(g(x)) =o(1)$.
	\end{lemma}
	\begin{remark}[Asymptotic Algebra]\label{R: Asymptotic Algebra} {\em
		It is customary in asymptotic analysis to write simply $o(x^{m})\pm o(x^{n})=o(x^{k})$, $o(x^{m})\pm o(x^{n})=o(x^{m+n})$ and $o(o(1)) = o(1)$ instead of (a), (b), and (c) in the above lemma, respectively, when no confusion could arise. We shall use this notation. 
	}\end{remark}
	\begin{lemma}[Absolute Value] \label{L: Norm} Let $A\in\mathbb{R}$. Then $|A+o(x^n)|=|A|+o(x^n)$. 
	\end{lemma}

\section{DEFINITE INTEGRALS WITHOUT RIEMANN SUMS}\label{S: Definite Integral Without Riemann's Sums}

	\begin{lemma}[Uniqueness]\label{L: Uniqueness}
		Let $a, b\in\mathbb{R},\; a<b$, and let $\rho: [a,b] \to \mathbb{R}$ be a continuous function.  Let $I: [a, b]\times[a, b]\to\mathbb{R}$ be a function in two variables satisfying the axioms {\bf (A)} and {\bf (B)} at the beginning of the introduction. Then $I(a, x)$ satisfies the initial value problem $\frac{d}{dx} I(a, x) = \rho(x)$,\; $I(a, a) = 0$ on the interval $(a, b)$. Consequently, there is no more than one function $I: [a, b]\times[a, b]\to\mathbb{R}$ satisfying the properties {\bf (A)} and {\bf (B)}.
	\end{lemma}
	\Proof Suppose that $I(x, y)$ satisfies \textbf{(A)} and \textbf{(B)}. We have 
	
	\begin{align}\notag
& \frac{d}{dx}I(a,x)= \lim_{h \to 0}\frac{I(a, x+h) - I(a, x)}{h} \overset{ \text{\textbf{(A)}} } {= }   \lim_{h \to 0} \frac{I(x, x+h)}{h} = \notag\\
&\lim_{h \to 0} \frac{I(x, x+h) - \rho(x)h + \rho(x)h}{h} =\lim_{h \to 0}\frac{I(x, x+h)-\rho(x)h}{h}  + \rho(x)  
  \overset{ \text{\textbf{(B)}} } {= }   0+\rho(x)=\rho(x).
	\end{align}
	
Also, \textbf{(A)} implies $I(a, a) + I(a, b) = I(a, b)$, so  $I(a, a) = 0$. Suppose that $J(x, y)$ is another function satisfying the axioms {\bf (A)} and {\bf (B)}, and let $\Delta(x, y)=I(x, y)-J(x, y)$. We have $\frac{d}{dx} \Delta(a, x) = 0$ and $\Delta(a, a) = 0$,  and hence $\Delta(a, x)=0$ for all $x$ in $(a, b)$. Consequently, for a given $\rho(x)$ there is no more than one $I(a, x)$. Next, \textbf{(A)} implies $I(x, y)=I(a, y)-I(a, x)$. Thus $I(a, x)$ uniquely determines $I: [a, b]\times[a, b]\to\mathbb{R}$. $\blacktriangle$
	
	Notice that  the uniqueness result presented above does not involve partitioning of the interval $[a, b]$ and Riemann sums. Rather, it is based on the more elementary result from calculus that ``every function with zero derivative on an interval is a constant.'' 

	The above lemma justifies the following definition:
	
\begin{definition}[Axiomatic Definition]\label{D: Axiomatic Definition} {\em 
		Let $\rho: [a,b] \to \mathbb{R}$ be a continuous function and let $I: [a, b] \times [a, b] \to \mathbb{R}$ be a function in two variables satisfying the axioms {\bf (A)} and {\bf (B)} at the beginning of the introduction. Then  the value $I(a, b)$ is called the integral of $\rho(x)$ over $[a, b]$. We shall use the usual notation  $\int_{a}^{b}\rho(x)dx\eqdef I(a, b)$.
}\end{definition}
	
	Note that axioms  \textbf{(A)} and  \textbf{(B)} in the beginning of the introduction are easily motivated and visualized (see Figure 1). We observe as well that if $a\leq\alpha\leq\beta\leq a$, then the restriction $I\upharpoonright[\alpha, \beta]\times [\alpha, \beta]$ also satisfies the axioms  {\bf (A)} and {\bf (B)}, and thus  $\int_{\alpha}^{\beta}\rho(x)dx= I(\alpha, \beta)$.
	
	\begin{theorem}[Fundamental Theorem]\label{T: Fundamental Theorem} 
		 As before, let $\rho(x)$ be a continuous function on $[a, b]$.
				
		{\bf (i)}  If the integral $I(a, x)$ of $\rho(t)$ over $[a, x]$ exists for every $x$ in $[a, b]$, then 
 $\frac{d}{dx} I(a, x)=\rho(x)$ on $(a, b)$.

		{\bf (ii)} If $R(x)$ is an anti-derivative of $\rho$ on $[a, b]$, then  $I(a, b)= R(b)-R(a)$.
	\end{theorem}	
	\Proof Part (i) follows directly from Lemma~\ref{L: Uniqueness}.
	For (ii), we have  $\frac{d}{dx} I(a, x) = \rho(x)$\textnormal{ by (i) and}$\frac{d}{dx} R(x) =\rho(x)$ by assumption. It follows that $\frac{d}{dx}[I(a, x) - R(x)] = 0$ on $(a, b)$, implying $R(x) = I(a, x) + C$ for some constant $C$. Hence, $R(b) - R(a) = [I(a, b) + C - (I(a, a) + C)]  = I(a, b)$, as required, since $I(a, a)=0$ by Lemma~\ref{L: Uniqueness}. $\blacktriangle$
	
	Notice that our simple theory assumes that the function $I(x, y)$, and thereby the integral $I(a, b)$, exists. The following is our first existence result (for the general existence result see Theorem~\ref{T: General Existence Result}).

	\begin{corollary}[Weak Existence]\label{C: Weak Existence} 
		If $\rho(x)$ has an antiderivative on $[a, b]$, then the integral of $\rho(x)$ over  $[a, b]$ exists. 
	\end{corollary}
	\Proof Let $R(x)$ be an antiderivative of $\rho(x)$. Then the function  $I: [a, b] \times [a, b] \to \mathbb{R}$ defined by $I(x, y)=R(y)-R(x)$ satisfies \textbf{(A)}  and  \textbf{(B)} (at the beginning of the introduction) and it is the only function which satisfies \textbf{(A)}  and  \textbf{(B)} by Lemma~\ref{L: Uniqueness}. The number $I(a, b)$ is the integral we are looking for. $\blacktriangle$
	\begin{remark}[Basic Properties of the Integral]{\em 
		The basic properties of the integral  follow immediately from part (ii) of Theorem~\ref{T: Fundamental Theorem} under the assumption that the integrals exist. For example, Theorem~\ref{T: Fundamental Theorem} implies the linear property $\int [c_1f(x)+c_2g(x)]dx = c_1\int f(x)dx + c_2\int g(x)dx$ provided that at least two of the three integrals  exist.
}\end{remark}

	Our simple yet rigorous theory presented so far is powerful enough to support most of the topics related to the integral and its applications in a typical beginning calculus course. To deal with integrals such as $\int_a^b e^{-x^2}dx$ and $\int_a^b \sin{(x^2)}dx$, we need a general existence result. This (and only this) is the place in our approach where we use partitions of the interval $[a, b]$ and Riemann sums. 

	Let $[a, b]$ be, as before, a closed interval in $\mathbb{R}$ and $x, y \in\mathbb{R},\, a
\leq x< y \leq b$. Recall that  a {\em partition} of $[x, y]$ is  a finite ordered set of the form $P=\{x_0, x_1,\dots,x_n\}$, where  $n\in\mathbb{N}$ and  $x=x_0<x_1<\dots<x_n=y$. We denote by $\mathcal{P}[x, y]$ the set of all partitions of $[x, y]$. Let $\rho: [a, b]\to\mathbb{R}$ be a continuous function and $P=\{x_0, x_1,\dots,x_n\}$ be a partition of $[x, y]$.  Recall that 
	\begin{align}
		&L(P)=\sum_{k=1}^n\; \left(\min_{x_{k-1}\leq t\leq x_k}{\rho(t)}\right)(x_k-x_{k-1}), \notag\\
		&U(P)=\sum_{k=1}^n\; \left(\max_{x_{k-1}\leq t\leq x_k}{\rho(t)}\right)(x_k-x_{k-1}),\notag
	\end{align}
are called the {\em lower and upper Darboux sums} of $\rho(t)$  determined by $P$, respectively. Let $x, y, z\in\mathbb{R},\, a \leq x< y<z \leq b$. Notice that if $P\in\mathcal{P}[x, y]$ and $Q\in\mathcal{P}[y, z]$, then $P\cup Q\in\mathcal{P}[x, z]$ and we have $L(P)+L(Q)=L(P\cup Q)$ and $U(P)+U(Q)=U(P\cup Q)$. The next result can be found in any contemporary textbook on Riemann integration.

\begin{lemma}\label{L: Partitions} Let $\rho: [a, b]\to\mathbb{R}$ be a continuous function. Then:

	{\bf (i)} For every two partitions $P$ and $Q$ of $[a,  b]$ we have 
	\[
	(\min_{[a, b]} \rho)(b-a)\leq L(P)\leq U(Q)\leq (\max_{[a, b]} \rho)(b-a).
	\]

	{\bf (ii)}  For every $\epsilon\in\mathbb{R}_+$ there exists a partition $P$ of $[a, b]$ such that $U(P)-L(P)<\epsilon$.
\end{lemma}

	\begin{theorem}[General Existence Result]\label{T: General Existence Result} 
		Let  $\rho: [a, b]\to\mathbb{R}$ be a continuous function. Then $\rho$ has both an integral and an antiderivative. 
	\end{theorem}
	
	\Proof Suppose, first, that $x, y\in\mathbb{R},\, a\leq x<y\leq b$. We observe that the set $\{L(P)\mid P\in\mathcal{P}[x, y]\}$ is bounded from above by the number $(\max_{x\leq t\leq y}{\rho(t)})(y-x)$. Thus 
$I(x, y)=\sup\{L(P)\mid P\in\mathcal{P}[x, y]\}$ exists in $\mathbb{R}$ by the completeness of $\mathbb{R}$. We intend to show that $I(x, y)$ satisfies the axioms \textbf{(A)}  and \textbf{(B)} at the beginning of the introduction.  We start with (\textbf{B}): we have $(\min_{x\leq t\leq y}{\rho(t)})(y-x)\leq I(x, y)$, by the definition of $I(x, y)$, since $P=\{x, y\}$ is a partition of the interval $[x, y]$. We let $y-x=h$ and the result is $(\min_{x\leq t\leq x+h}{\rho(t)})h\leq I(x, x+h) \leq(\max_{x\leq t\leq x+h}{\rho(t)})h$. The latter implies axiom \textbf{(B)} since $\rho$ is continuous at $x$. To prove (\textbf{A}), we observe that  $L(P) \leq I(x, y)\leq U(P) $  for every partition $P\in\mathcal{P}[x, y]$. Indeed, the first inequality follows directly from the definition of $I(x, y)$ and the second inequality follows from the definition of $I(x, y)$ and part~(i) of Lemma~\ref{L: Partitions}. Next, suppose that $x, y, z\in\mathbb{R},\, a\leq x<y<z\leq b$, and let  $Q\in\mathcal{P}[y, z] $. As before we have $L(Q) \leq I(y, z)\leq U(Q) $, and also $L(P\cup Q) \leq I(x, z)\leq U(P\cup Q)$ since $P\cup Q$ is a partition of $[x, z]$. After summing up we obtain:
\[
L(P)+L(Q)-U(P\cup Q) \leq I(x, y)+I(y, z)-I(x, z)\leq U(P)+U(Q)-L(P\cup Q).
\]
Now we can choose partitions $P\in\mathcal{P}[x, y]$ and $Q\in\mathcal{P}[y, z] $ such that  $U(P)-L(P)<\epsilon/2$ and $U(Q)-L(Q)<\epsilon/2$ by part~(ii) of Lemma~\ref{L: Partitions}. Also, $U(P\cup Q)-L(P\cup Q)<\epsilon$ since $L(P)+L(Q)=L(P\cup Q)$ and $U(P)+U(Q)=U(P\cup Q)$.  Thus the above chain of inequalities reduces to 
$-\epsilon<I(x, y)+I(y, z)-I(x, z)<\epsilon$, implying $I(x, y)+I(y, z)=I(x, z)$, as required. Finally we can eliminate the restriction on $x, y$ and $z$ by extending the function $I(x, y)$ to a function in the form $I: [a, b]\times [a, b]\to\mathbb{R}$ by letting $I(x, y)=-I(y, x)$ and $I(x, x)=0$ for all $x$ and $y$ in $[a, b]$. Notice that the function $I(x, y)$ just defined also satisfies {\bf (A)} and {\bf (B)};  thus it is uniquely determined, by  Lemma~\ref{L: Uniqueness}. The number $I(a, b)$ is the integral we are looking for, by Definition~\ref{D: Axiomatic Definition}, and $R(x)=I(a, x)$ is an antiderivative of $\rho(x)$, by part ~(i) of Theorem~\ref{T: Fundamental Theorem}. $\blacktriangle$

\section{SETTING UP INTEGRALS WITHOUT RIEMANN SUMS}\label{S: Setting Up Integrals Without Riemann's Sums}

	Imagine that you are a young instructor preparing to cover  arclength in a typical calculus course. We can safely assume that before going in front of your students you would like to clarify the structure of the topic for yourself: what is the definition of ``arclength,''  how do I motivate it, is there a theorem to present (with or without formal proof), and finally which examples should I use? One option is to define the concept of arclength directly by the integral formula $L(a, b) = \int_{a}^{b} \sqrt{1 + f'^{2}(x)}dx$. If this is your choice, your next task will be to motivate this definition. You might use the integral formula to calculate some familiar results from high school mathematics: the distance formula  for a line segment or the formula for the circumference of a circle.  This approach, although completely legitimate, is rarely used by calculus textbooks. The integral formula still looks terribly unmotivated even after deriving the distance formula for a line segment. Also, it is far from clear that this is the only formula producing the distance formula.  For that reason most  calculus books use Riemann sums to convince students that the integral formula is ``reasonable.''  What follows is messy mathematics: the partition of the interval $[a, b]$ gives the impression that the step of the partition $\Delta x$  is ``very small''  which leads to the conclusion that $\Delta L$ is approximately equal to $\sqrt{\Delta x^2 +\Delta y^2}$. The latter produces our integral formula after factoring out $\Delta x$ and taking the limit as $\Delta x\to 0$. You might ponder for hours questions such as: a) What, after all, is the definition of ``arclength''? Is the integral formula {\em exact or approximate}?  After all, its derivation involves the approximate formula $\Delta L\approx\sqrt{\Delta x^2 +\Delta y^2}$. b) What is meant by ``$\Delta x$ is very small"? If $\Delta x=0$, then $\Delta y=0$ and the root $\sqrt{\Delta x^2 +\Delta y^2}$ is also equal to zero. If $\Delta x\not=0$, then  $\Delta L\not=\sqrt{\Delta x^2 +\Delta y^2}$ (unless the curve is a straight line).  c) Is the ``derivation" of the formula for arclength a sort of casually presented proof of a casually stated theorem? And if yes, what is the rigorous version of this theorem? Worst of all is that the students are usually preoccupied with the technicalities of the partition procedure and the sigma notation in the Riemann sums and hardly pay attention to the fact that a new important concept has just been introduced.
	
	In this section we take another approach.  The concept of {\em arclength} is defined as an additive quantity which is asymptotically equal to the Euclidean distance between two points (Definition~\ref{D: Arc-Length}). The definition is mathematically correct and well motivated. It is based on the concept of limit, not integral. Next we derive the arclength integral formula as a simple rigorous theorem in the spirit of asymptotic analysis (Theorem~\ref{T: Arc-Length}). Similarly we define the rest of the additive quantities from geometry and physics such as  {\em area under a curve,  volume of a solid of revolution}, etc. We involve neither a partition of the interval $[a, b]$ nor Riemann sums.
	
	We start with area under a curve. Our assumption is that the reader knows what the {\em area of a rectangle} is but not necessarily what the  {\em area under a curve} is. In particular, we do not assume that the reader necessarily  knows the integral formula $A(a, b) = \int_{a}^{b} f(x)d\,x$ for the area under a curve; rather our goal is to derive this formula starting from the more elementary concept of the area of a rectangle.

	\begin{definition}[Area Under a Curve]\label{D: Area}{\em 
		Let $y = f(x)$ be a continuous function on $[a, b]$ such that $f(x) \geq 0$ for all $x\in [a, b]$. Let  $A: [a, b] \times [a, b] \to \mathbb{R}$ be a function in two variables satisfying the following two properties: 
		\begin{description}
			\item{\bf (a)} $A(x, y) + A(y, z) = A(x, z)$ for all $x, y, z \in [a, b]$.
			\item{\bf (b)} $A(x, x + h) = \pm R(x, x + h) + o(h)$ as $h\to 0_\pm$ for all $x \in [a, b]$, where $R(x, x+h)$ denotes the area of the rectangle with vertices $(x, 0), (x+h, 0),
(x+h, f(x))$, and $(x, f(x))$.
		\end{description}
The number $A(a, b)$ is called (by definition)  the {\bf area under the curve} $y=f(x)$ and  above the interval $[a, b]$.
	}\end{definition}
	
		The above definition can be easily motivated (see Figure 2).
	
	
	In the next theorem we derive the familiar integral formula for $A(a, b)$ without partitions or Riemann sums. While deriving this formula, we shall simultaneously prove two things: (a)  the correctness of the above definition, and (b)  the existence of the area under the curve.  As in the conventional approach, the integral formula offers a practical method for explicit evaluation. 	
	\begin{theorem}
		$A(a, b) = \int_{a}^{b} f(x)dx$. Consequently, the area $A(a, b)$ is uniquely determined by the properties  (a) and (b) in Definition ~\ref{D: Area}.
	\end{theorem}
	\Proof We have to find the asymptotic expansion of $A(x, x+h) $ in powers of $h$ as $h\to 0$ and extract the coefficient in front of $h$. Since $R(x, x+h) = f(x)|h|$, we have $A(x, x+h) = \pm R(x, x+h) + o(h) =\pm f(x)|h| + o(h)= f(x)h + o(h)$, and the above formula follows directly from Definition~\ref{D: Axiomatic Definition} for $\rho(x) = f(x)$. Conversely, it is easy to verify that the function  $A: [a, b] \times [a, b] \to \mathbb{R}$ defined by  $A(x,y) = \int_{x}^{y} f(t)d\,t$ satisfies the properties (a) and (b) in the above definition, thus the number $A(a, b)$ is the area under the curve. $\blacktriangle$
	
	Next,  we define the concept of {\em arclength} without partitions or Riemann sums.  Our assumption is that the reader knows what the {\em Euclidean distance between two points} is but not necessarily what the  {\em arclength along a curve} is. In particular, we do not assume any knowledge about the integral formula $L(a, b) = \int_{a}^{b} \sqrt{1 + f'^{2}(x)}dx$; our goal is to derive this formula starting from the more elementary concept of  Euclidean distance between two points.
	
	\begin{definition}[Arclength] \label{D: Arc-Length} {\em  Let $f \in C^{1}[a, b]$ and let $L: [a, b] \times [a, b]\to \mathbb{R}$ be a function in two variables satisfying the following two properties:
		\begin{description}
			\item{\bf (a)} $L(x, y) + L(y, z) = L(x, z)$ for all $x, y, z \in [a, b]$.
			\item{\bf (b)} $L(x, x + h) = \pm D(x, x + h) + o(h)$ as $h\to 0_\pm$ for all $x \in [a, b]$, where $D(x, x+h)$ is the Euclidean distance between the points $(x, f(x))$ and $(x+h, f(x+h))$.
		\end{description}
 The number $L(a, b)$  is called (by definition) the {\bf arclength} of the curve $y=f(x)$ between the points $(a, f(a))$ and $(b, f(b))$.
	}\end{definition}
	
	The above definition can be easily motivated (see Figure 3).


	 In the next theorem we rigorously derive the formula for  arclength without involving partitions of the interval or Riemann sums. Among other things we prove correctness of the above definition and the existence of the arclength $L(a, b)$. 
	\begin{theorem}\label{T: Arc-Length}
		$L(a, b) = \int_{a}^{b} \sqrt{1 + f'^{2}(x)}dx$. Consequently, the arc length $L(a, b)$ is uniquely determined by the properties  (a) and (b) in Definition~\ref{D: Arc-Length}.
	\end{theorem}		
	\Proof  We have to find the asymptotic expansion of  $L(x, x + h)$ in powers of $h$ as $h\to 0$ and extract the coefficient $\rho(x)$ in front of $h$ (see Definition~\ref{D: Axiomatic Definition}). Let $\Delta y = f(x + h) - f(x)$ and recall that $\Delta y = f'(x)h + o(h)$ (Increment Theorem~\ref{T: Increment Theorem}). We have: 
		\begin{align}\notag
&L(x, x + h)=\pm D(x, x+h)+o(h) = \pm \sqrt{h^2 + {\Delta y}^2} + o(h)=\\\notag
&=\pm |h| \sqrt{1 + \left(\frac{f^{\prime}(x)h + o(h)}{h}\right)^{2}} +o(h)=\\\notag
&=h\sqrt{1+f^{\prime 2}(x)+\frac{2f^{\prime}(x)o(h)h}{h^2}
+\left(\frac{o(h)}{h}\right)^2}+o(h)=\\\notag   
&=h\sqrt{1 +f^{\prime 2}(x) + o(1)} + o(h)=h\left[\sqrt{1+f^{\prime 2}(x)}+\sqrt{1 + f^{\prime 2}(x)
+o(1)}-\sqrt{1+f^{\prime 2}(x)} \right] + \\\notag
&+o(h)=h\left[\sqrt{1+f^{\prime 2}(x)}+\frac{o(1)}{\sqrt{1+f^{\prime
2}(x)+o(1)}+\sqrt{1+f^{\prime 2}(x)}}\right]+o(h)=\\\notag   
&=h\left[\sqrt{1+f^{\prime
2}(x)}+o(1)\right]+o(h)=h\sqrt{1+f^{\prime 2}(x)} +o(h) + o(h)=\\\notag  
&=h\sqrt{1+f^{\prime 2}(x)} + o(h)\notag.
\end{align}  
	For the coefficient in front of $h$ we have $\rho(x) = \sqrt{1 + f'^{2}(x)}$, which implies our integral formula by Definition~\ref{D: Axiomatic Definition}.   Conversely, it is easy to verify that the  function $L: [a, b] \times [a, b]\to \mathbb{R}$ defined by $L(x, y) = \int_{x}^{y} \sqrt{1 + f'^{2}(t)}d\, t$ satisfies (a) and (b) in Definition~\ref {D: Arc-Length} . Thus the number  $L(a, b)$ is the arc length of the curve.   $\blacktriangle$

		Next, we set up the integral for a {\em volume of revolution} about the $y$-axis without partitions of the interval and Riemann sums. Our assumption is that the reader knows what a  {\em volume of a cylindrical shell} is but not necessarily what a  {\em volume of revolution} is.
	
	\begin{definition}[Volume of Revolution] \label{D: Shell} {\em 
		Let $y = f(x)$ be continuous on $[a, b]$, $f(x) \geq 0$, and $0 \leq a < b$. Let $V: [a, b] \times [a, b] \to \mathbb{R}$ be a function in two variables satisfying the following properties:
			
		{\bf (a)} $V(x, y) + V(y, z) = V(x, z)$ for all $x, y, z \in [a, b]$.
			
		{\bf (b)} $V(x, x + h) = \pm U(x, x+h) + o(h)$ as $h \to 0_{\pm}$ for all $x\in [a, b]$, where $U(x, x+h)$ is the volume of the cylindrical shell obtained by revolving the rectangle with vertices $(x, 0),  (x + h, 0), (x+h, f(x))$ and $(x, f(x))$ about the $y$-axis (see Figure 4).
		
	The number $V(a, b)$ is called (by definition) the {\bf volume of revolution about the $y$-axis} of the curve $y=f(x)$.
	}\end{definition}
	
	 In the next theorem we rigorously derive the familiar formula for the volume $V(a, b)$, and we do so without partitions or Riemann sums. Among other things we prove the uniqueness and existence of the volume $V(a, b)$.	
	\begin{theorem}
		$V(a, b) = \int_{a}^{b} 2\pi x f(x)dx$. Consequently, the volume $V(a, b)$ is uniquely determined by the properties  (a) and (b) from Definition ~\ref{D: Shell}.
	\end{theorem}
	\Proof The volume of the cylindrical shell is $U(x, x+h) = |\pi (x + h)^{2} - \pi x^{2}|f(x)$. Hence, with the help of Lemma~\ref{L: Norm}, we have $V(x, x + h) = \pm U(x, x+h) +o(h) = \pm |\pi (x + h)^{2} - \pi x^{2}|f(x) + o(h) = \pm |\pi (2xh + h^{2})|f(x) + o(h) = \pm 2\pi xf(x)|h| + \pi f(x) h^{2} + o(h) = \pm 2\pi xf(x)|h| + o(h) + o(h) = 2\pi xf(x)h + o(h)$, and from Definition~\ref{D: Axiomatic Definition} we have $\rho(x) = 2\pi xf(x)$.  Conversely, it is easy to verify that the function  $V: [a, b] \times [a, b] \to \mathbb{R}$ defined by   $V(x, y) = \int_{x}^{y} 2\pi t f(t)dt$ satisfies the properties (a) and (b) in Definition~\ref{D: Shell}. Thus the number $V(a, b)$ is the volume of revolution. $\blacktriangle$
	

\section{BACK TO INFINITESIMAL CALCULUS: A LOST THEOREM} \label{S: Lost Theorem}

	As we explained at the end of the introduction, the method presented in this article has a long and interesting history. The purpose of this section is to establish a relation between our method of integration and infinitesimal calculus. This section may be helpful to those readers who are interested in reading original texts on infinitesimal calculus but who might not have background in modern nonstandard analysis. We should mention that the tradition of using infinitesimal arguments  is still very much alive and can be found in many contemporary texts on applied mathematics, physics, and engineering. For example, in the famous {\em The Feynman Lectures on Physics}~\cite{Feynman} we located  about ten cases of integral formulas for additive physical quantities derived in the spirit of infinitesimal calculus and without Riemann sums (see Volume I: pp. 13-3, 14-8, 43-2, 44-10 and 44-11, 46-6, 47-5; Volume II: pp. 3-2, 38-6). For that reason we believe that our article and the discussion in this section in particular might also be helpful to pure mathematicians who are interested in reading texts on applied mathematics, physics, or engineering but who might feel uneasy with infinitesimal  reasoning.

	 Recall that in the period from Leibniz to Weierstrass,  calculus was commonly known as  {\em infinitesimal calculus} and was based on the hypothesis that there exist nonzero infinitesimals, i.e., mysterious numbers $dx$ with  the property $0<|dx|<1/n$ for all $n\in\mathbb{N}$.  We should keep in mind that at that period not only the theory of infinitesimals but also the theory of real numbers was without rigorous foundation. So, the existence of nonzero infinitesimals should not be dismissed as nonsense.  Yes, the field of the real numbers $\mathbb{R}$ does not have nonzero infinitesimals,  but there were no {\em real numbers} in the era of Leibniz and Euler;  the real numbers were an invention of the late 19th century and were systematically implemented in mathematics at the beginning of 20th century. 

	We shall discuss Leibniz-Euler infinitesimal calculus using the common ``differential notation'': Let $y=f(x), a\leq x\leq b$ be a function. In what follows, $dx$  stands for a new independent variable (real or infinitesimal depending on the context) and $\Delta y=f(x+dx)-f(x)$ stands for the increment of $y$. If $f$ is differentiable at $x$, then $dy=f^\prime(x)dx$ stands for the differential of $y$. It is clear that  $dx=\Delta x$, but we prefer to use $dx$ over $\Delta x$, thus keeping track of the fact that  $x$ is an independent and $y$ is a dependent variable. Before proceeding further we should notice that in the old infinitesimal calculus, and in the modern nonstandard analysis as well, the notation $dx$ rarely stands for a {\em fixed infinitesimal number} (as $\pi, e, \sqrt{2}$, etc. stand for specific irrational numbers). Rather, $dx$ is usually used for an independent variable ranging over a set consisting of  both infinitesimal and real (standard) numbers. To demonstrate this point we shall present the characterization of  continuity used by Euler but rigorously justified in the modern nonstandard analysis:  Let $x$ be a real (standard) number in the domain of $f$. Then $f$ is continuous at $x$ if and only if $f(x+dx)-f(x)$ is infinitesimal for every infinitesimal $dx$. In the manuscripts of Euler this statement appears in a slightly more casual form: $f$ is continuous at $x$ if 
and only if $\Delta y$ is infinitesimal whenever $dx$ is infinitesimal. It is clear that if both $x$ and $dx$ are real (standard) numbers, then $\Delta y=f(x+dx)-f(x)$ is also a real (standard) number. Also, ``real (standard) number'' is a modern term. Leibniz and Euler would rather use ``usual quantity'' (as opposed to ``infinitesimal quantity'') instead.

	 Let us try to mimic, for example, the arguments used by Euler for setting up the integral for the arclength $L(a, b)$ along the curve $y=f(x), a\leq x\leq b$ in a typical calculus course in the middle of the 18th century (compare with our Definition~\ref{D:  Arc-Length} and Theorem~\ref{T: Arc-Length}). If $dx$ is infinitesimal, the arclength $L(x, x+dx)$ between the points $(x, f(x))$ and $(x+dx, f(x+dx))$ should be {\em equal} to the Euclidean distance $D(x, x+dx)=\sqrt{dx^2+\Delta y^2}$ between the same points, {\em up to infinitesimals of second order relative to $dx$};  in symbols, $L(x, x+dx)\approx \sqrt{dx^2+\Delta y^2}$ (see Figure 3, where $h$ should be replaced by $dx$).  On the other hand, since $L$ is an additive quantity, we have $L(x, x+dx)=L(a, x+dx)-L(a, x)$. Thus $L(x, x+dx)\approx dL$, by the increment theorem (in its original infinitesimal form).  Also, $\Delta y\approx dy$, by the increment theorem. The result is $dL\approx \sqrt{1+(dy/dx)^2}\, dx$ which implies the familiar integral formula $L(a,b)=\int_a^b\sqrt{1+(dy/dx)^2}\, dx$. We should note that the implication in the previous sentence has never been  rigorously justified in the old infinitesimal calculus; among other things the goal of our article is to fill this gap. 
	 	 
	A contemporary mathematician, unless familiar with nonstandard analysis, will have difficulty recognizing our asymptotic property {\bf (B)} (at the beginning of the introduction) using Euler's language presented above. For those who are interested in using nonstandard analysis we recommend   H. J. Keisler~\cite{Keisler} or T. Todorov~\cite{tdTod2000a}, where the reader will find more literature on the subject. In what follows, however, we shall choose another path: we shall use the language of asymptotic analysis (Section~\ref{S: Topics in Asymptotic Analysis}) to relate the method of integration presented here with the method of infinitesimal calculus. For that purpose we suggest the following modification of Section~\ref{S: Definite Integral Without Riemann's Sums} and Section~\ref{S: Setting Up Integrals Without Riemann's Sums} of our article: 

	1) First, instead of the letter $h$ (used in our article so far) we shall use the original Leibniz notation $dx$. In other words we let $dx=h$ and we treat $dx$ as a new real (standard) independent variable. In this notation our Asymptotic Property {\bf (B)} (see the beginning of the introduction) becomes
	\begin{equation}\notag
		{\bf (B^\prime)\; \;} I(x, x + dx) = f(x)dx + o(dx)  {\; \text  as\; } dx\to 0  \text{\; for all\;} x \in [a, b].
	\end{equation}
As before, $f(dx)=o(dx^n)$ means $f(dx)/dx^n\to 0$ as $dx\to 0$. 

	2) Our next goal is to give precise meaning to the relation $\approx$ used by Euler in our previous example. Let $F(dx)$ and $G(dx)$ be two real functions. We say that $F(dx)$ and $G(dx)$ are {\em equal up to infinitesimals of second order relative to $dx$}, in symbols, $F\approx G$, if $F(dx) -G(dx)=o(dx)$ (i.e., if $ (F(dx) -G(dx))/dx\to 0$ as $dx\to 0$). The property  {\bf (B$^\prime$)} becomes
	\begin{equation}\notag
		{\bf (B^{\prime\prime})\; \;} I(x, x + dx) \approx f(x)dx \text{\; for all\;} x \in[a, b].
	\end{equation}
	
	3) We might stop writing  $dx\to 0$, since the symbol $dx$ is more than suggestive.  
		
	4) Now we can rewrite Section~\ref{S: Definite Integral Without Riemann's Sums},  replacing $h$ by $dx$ and the axiom {\bf (B)} by  $(\mathbf{B}^{\prime\prime})$. The additive property {\bf (A)}  (at the beginning of the introduction) does not need modification. 
	
	5) Finally, we have to replace all parts~{\bf (b)} in the definitions in Section~\ref{S: Setting Up Integrals Without Riemann's Sums} by their ({\bf b}$^{\prime\prime}$)-counterparts in the spirit of $(\mathbf{B}^{\prime\prime})$. For example, part~{\bf (b)} of Definition~\ref{D: Arc-Length} should  be replaced by:
\[
	{\bf (b^{\prime\prime})} \quad  L(x, x + dx) \approx \pm D(x, x + dx) \text{\;for all\; } x \in [a, b],
\]
and similarly for the rest of the ({\bf b)s} in Section~\ref{S: Setting Up Integrals Without Riemann's Sums}.

	While preserving the content of the article, the modification suggested above makes the 
arguments of the old infinitesimal calculus more transparent.

	6) Readers who feel uncomfortable with the notation $dx$ (and especially with $dx\to 0$) should restore the notation $h$ used in Section 3 and Section 4 of this article. However we recommend our trick ``replace $h$ by $dx$'' to those readers who are interested in reading original texts on infinitesimal calculus, but do not have background in the modern nonstandard analysis.
	
	 Among other things our article suggests that the reasoning of mathematicians  in the era of  Leibniz  and  Euler, as well as  the reasoning of  contemporary applied mathematicians and  physicists, is often more  reliable and rigorous than is usually acknowledged by pure mathematicians.\\

{\bf Acknowledgement:} We are thankful to our colleague Donald Hartig who made several useful remarks on the manuscript.
}

\begin{description}
\item[Ray Cavalcante] was an undergraduate student  in mathematics at the time this article was written, and he is currently a graduate student at Cal Poly, San Luis Obispo. 

 Mathematics Department, California Polytechnic State University, San Luis Obispo, CA 93407, USA (raycavalcante@gmail.com).
\item[Todor D. Todorov]  received his Ph.D. in Mathematical Physics from University of Sofia, Bulgaria. He currently teaches mathematics at Cal Poly, San Luis Obispo. His articles are on nonstandard analysis, nonlinear theory of generalized functions (J. F. Colombeau's theory), asymptotic analysis, compactifications of ordered topological spaces, and teaching calculus. Presently he works on a nonstandard version of the nonlinear theory of generalized functions and its applications to PDE. 

 Mathematics Department, California Polytechnic State University, San Luis Obispo, CA 93407, USA  (ttodotrov@calpoly.edu).
\end{description}
\end{document}